\documentclass[10pt,a4paper,leqno]{amsart}
\usepackage[english]{babel}
\usepackage[T1]{fontenc}
\usepackage{amssymb}
\usepackage{amsfonts}
\usepackage{amsmath}
\usepackage{mathptm}

\usepackage{eucal}
\begin{document}

%
%

\title{\large{ On Separately Subharmonic and Harmonic Functions}}
\author{Juhani Riihentaus}
\date{February 11, 2013}
\maketitle
%
\begin{center}
{Department of Mathematical Sciences, University of Oulu\\
P.O. Box 3000, FI-90014 Oulun yliopisto, Finland \\
riihentaus@member.ams.org}
\end{center}
\begin{center}
{and}
\end{center}
\begin{center}
{Department of Physics and Mathematics, University of Eastern Finland\\
P.O. Box 111, FI-80101 Joensuu, Finland \\
juhani.riihentaus@uef.fi}
\end{center}

\vspace*{2ex}

\vspace*{2ex}

\noindent{\emph{Abstract:}} We improve our previous generalizations to Arsove's and Ko\l odziej's and Thorbi\"ornson's results concerning the subharmonicity of a function subharmonic with respect to the first variable and harmonic with respect to the second.

\vspace{0.3ex}


\vspace{0.3ex}

\noindent{{{\emph{AMS 2010  Mathematics Subject Classification:} 31C05, 31B25, 31B05
\section{Introduction}
\subsection{ } Solving a long standing problem, Wiegerinck [Wi88, Theorem, p.~770], see also Wiegerinck and Zeinstra [WiZe91, Theorem~1, p.~246], showed that a separately subharmonic function need not be subharmonic. On the other hand, it is an open problem, whether a function which is subharmonic in
 one variable and harmonic in the other, is subharmonic. For older results on this area,
see e.g. Arsove [Ar66, Theorem~2, p.~622], Imomkulov [Im90, Theorem, p.~9], Wiegerinck and Zeinstra [WiZe91, p.~248], Cegrell and Sadullaev [CeSa93, Theorem~3.1, p.~82] and Ko\l odziej and Thorbi\"ornson [KoTh96, Theorem~1, p.~463].   The result of Ko\l odziej and Thorbi\"ornson includes the results of Arsove, of Cegrell and Sadullaev and of Imomkulov, and  reads as  follows:

\vspace{0.3ex}

\noindent{\textbf{Theorem~A}} \emph{ Let $\Omega $ be a domain in ${\mathbb{R}}^{m+n}$, \mbox{$m,n\geq 2$.}
Let $u:\, \Omega \rightarrow
{\mathbb{R}}$ be such that}
 \begin{itemize}
\item[(a)] \emph{for each $y\in {\mathbb{R}}^n$ the function}
\[\Omega (y)\ni x\mapsto u(x,y)\in {\mathbb{R}}\]
\emph{is subharmonic and ${\mathcal{C}}^2$,}
 \item[(b)] \emph{for each $x\in {\mathbb{R}}^m$ the function}
\[\Omega (x)\ni y\mapsto u(x,y)\in {\mathbb{R}}\]
\emph{is harmonic.}
\end{itemize}
\emph{Then} $u$
\emph{is subharmonic and continuous in $\Omega $}.

\vspace{0.3ex}

We improved the result of Ko\l odziej and Thorbi\"ornson in a series of papers: [Ri07$_1$, Theorem~3, Theorem~4 and Corollary, pp.~162--164], [Ri07$_2$, Theorem~6, p. 234],  [Ri07$_3$, Theorem~1 and Corollary, pp.~438, 444], [Ri07$_4$, Theorem~5.1,  Corollary~5.1 and Corollary~5.2, pp.~67--68, 74] and [Ri09, Theorem~4.3.1,   Corollary~4.3.3 and Corollary~4.3.4, pp.~e2625--e2626].  We will now return to the subject and improve our result still further, see Theorem~2 below. 
\subsection{ } However, we begin with improving the above cited results of Arsove and of Cegrell and Sadullaev  and  our previous generalizations  [Ri07$_4$, Theorem~4.1, p.~64] and [Ri09, Theorem~4.2.1, p.~e2623]. Instead of subharmonic functions (resp. so called quasinearly subharmonic functions n.s.),  we will now use quasinearly subharmonic functions.  Observe that in certain situations such an approach is indeed useful. One such an example is the following. Armitage and Gardiner [ArGa93, Theorem~1, p.~256] gave a condition which ensures a separately subharmonic function to be subharmonic,  and this  condition  was close to being sharp, see [ArGa93, pp.~255--256]. With the aid of quasinearly subharmonic functions it was, nevertheless, possible to generalize and improve their result,  see [Ri08, Theorem~4.1 and Corollary~4.5, pp.~8-9,~13] and [Ri09, Theorem~3.3.1 and Corollary~3.3.3,  pp.~e2621--e2622].   
\subsection{ } Our presentation below, including the presented references, is rather detailed. For the notation, and for the definitions and properties of subharmonic functions,  nearly subharmonic functions,  quasinearly subharmonic functions (and quasinearly subharmonic functions n.s., too) etc., see e.g. [Br69], [He71], [Ta88], [RiTa93], [Ri07$_4$], [Ri08], [PaRi08], [Ri09], [PaRi09], [Ri11],  and the references therein. 
\section{Arsove's result and its improvement}  
\subsection{ }  Arsove's result is:

\vspace{0.3ex}

\noindent{\textbf{Theorem~B}} ([Ar66, Theorem~2, p.~622]) \emph{ Let $\Omega $ be a domain in ${\mathbb{R}}^{m+n}$, \mbox{$m,n\geq 2$.}
Let $u:\, \Omega \rightarrow
{\mathbb{R}}$ be such that}
 \begin{itemize}
\item[(a)] \emph{for each $y\in {\mathbb{R}}^n$ the function}
\[\Omega (y)\ni x\mapsto u(x,y)\in {\mathbb{R}}\]
\emph{is subharmonic,}
 \item[(b)] \emph{for each $x\in {\mathbb{R}}^m$ the function}
\[\Omega (x)\ni y\mapsto u(x,y)\in {\mathbb{R}}\]
\emph{is harmonic,}
\item[(c)] \emph{there is a nonnegative function } $\varphi \in {\mathcal{L}}^1_{\textrm{loc}}(\Omega )$ \emph{such that} $-\varphi \leq u$.
\end{itemize}
\emph{Then} $u$
\emph{is subharmonic in $\Omega $}.

\vspace{0.3ex}

Arsove's proof was based on mean value operators. Much later  Cegrell and Sadullaev
[CeSa93, Theorem~3.1, p. 82] gave a new proof using Poisson modification.
\subsection{ } Below in Theorem~1 we generalize the above result of Arsove and of  Cegrell and Sadullaev, and also our previous generalizations  [Ri07$_4$, Theorem~4.1 and Corollary~4.1, pp.~64--65], see also [Ri09, Theorem~4.2.1 and Corollary~4.2.2, p.~e2623]. The proof we give below is  short and direct. 

\vspace{0.3ex}

\noindent{\textbf{Theorem~1}} \emph{Let $\Omega $ be a domain in ${\mathbb{R}}^{m+n}$, \mbox{$m,n\geq 2$, and $K\geq 1$.}
Let $u:\, \Omega \rightarrow
{\mathbb{R}}$ be  such that}
 \begin{itemize}
\item[(a)] \emph{for each $y\in {\mathbb{R}}^n$ the function}
\[\Omega (y)\ni x\mapsto u(x,y)\in {\mathbb{R}}\]
\emph{is K-quasinearly subharmonic,}
 \item[(b)] \emph{for each $x\in {\mathbb{R}}^m$ the function}
\[\Omega (x)\ni y\mapsto u(x,y)\in {\mathbb{R}}\]
\emph{is harmonic,}
\item[(c)] \emph{there is a nonnegative function } $\varphi \in {\mathcal{L}}^1_{\textrm{loc}}(\Omega )$ \emph{such that} $-\varphi \leq u$.
\end{itemize}
\emph{Then} $u$
\emph{is K-quasinearly subharmonic in $\Omega $}.

\vspace{0.3ex}

\noindent\emph{Proof.} It is easy to see that $u$ is Lebesgue measurable.  Therefore also $u_M:=\max \{\,u,-M\,\}+M$, $M>0$, is Lebesgue measurable. We must show that $u^+\in {\mathcal{L}}^1_{\textrm{loc}}(\Omega)$ and that each $u_M$ satisfies the generalized mean value inequality. 

To see that $u^+\in {\mathcal{L}}^1_{\textrm{loc}}(\Omega)$, we proceed as follows. Observe first that $0\leq u^+\leq u_M\leq v_M:=u+\varphi +M$. To see that $v_M\in {\mathcal{L}}^1_{\textrm{loc}}(\Omega)$ requires \emph{only} Fubini's Theorem. As a matter of fact, take  $\overline{B^m(a,R)\times B^n(b,R)}\subset \Omega $ arbitrarily.
Then
\begin{displaymath}\begin{split}0\leq& \frac{K}{m_{m+n}(B^m(a,R)\times B^n(b,R))}\int\limits_{B^m(a,R)\times B^n(b,R)}v_M(x,y) dm_{m+n}(x,y)=\\
 \leq& \frac{K}{m_{m+n}(B^m(a,R)\times B^n(b,R))}\int\limits_{B^m(a,R)\times B^n(b,R)}[u(x,y)+\varphi (x,y)+M] dm_{m+n}(x,y)=\\
  \leq&\frac{K}{\nu _m\, R^m}\int\limits_{B^m(a,R)}\{\frac{1}{\nu _n\,R^n}\int\limits_{B^n(b,R)} [u(x,y)+\varphi (x,y)+M] dm_n(y)\}dm_m(x)=\\
\leq&\frac{K}{\nu _m\, R^m}\int\limits_{B^m(a,R)}[\frac{1}{\nu _n\,R^n}\int\limits_{B^n(b,R)} u(x,y)dm_n(y)+
\frac{1}{\nu _n\,R^n}\int\limits_{B^n(b,R)} \varphi (x,y) dm_n(y)+M]dm_m(x)=\\
\leq&\frac{K}{\nu _m\, R^m}\int\limits_{B^m(a,R)}[u(x,b)+
\frac{1}{\nu _n\,R^n}\int\limits_{B^n(b,R)} \varphi (x,y) dm_n(y)+M]dm_m(x)=\\
\leq&\frac{K}{\nu _m\, R^m}\int\limits_{B^m(a,R)}u(x,b)dm_m(x)+\frac{K}{\nu _mR^m}\int\limits_{B^m(a,R)}[
\frac{1}{\nu _n\,R^n}\int\limits_{B^n(b,R)} \varphi (x,y) dm_n(y)]dm_m(x)+K\, M=\\
\leq&\frac{K}{\nu _m\, R^m}\int\limits_{B^m(a,R)}u(x,b)dm_m(x)+\frac{K}{m_{m+n}(B^m(a,R)\times B^n(b,R))}\int\limits_{B^m(a,R)\times B^n(b,R)}
\varphi (x,y) dm_{m+n}(x,y)+K\, M\\
<&+\infty .
\end{split}\end{displaymath}

It remains to show that for all  $(a,b)\in \Omega $ and $R>0$  such that $\overline{B^{m+n}((a,b),R)}\subset \Omega $,
\begin{displaymath}u_M(a,b)\leq \frac{K}{\nu _{m+n}R^{m+n}}\int\limits_{B^{m+n}((a,b),R)}u_M(x,y)dm_{m+n}(x,y).\end{displaymath}
To see this, we proceed in the following standard, direct and short way, see e.g. [He71, Proposition~2~c) and proof of Theorem~a), pp.~10--11, 32--33] and [Ri07$_4$, p.~59]:
\begin{align*}&\frac{K}{\nu _{m+n}R^{m+n}}\int\limits_{B^{m+n}((a,b),R)}u_M(x,y)dm_{m+n}(x,y)=\\
&=\frac{\nu _m}{\nu _{m+n}R^{m+n}}\int\limits_{B^n(b,R)}[(R^2-\mid y-b\mid ^2)^{\frac{m}{2}}\frac{K}{\nu _m(R^2-\mid y-b\mid ^2)^{\frac{m}{2}}}
\int\limits_{B^m(a,\sqrt{R^2-\mid y-b\mid ^2})}u_M(x,y)dm_m(x)]dm_n(y) \\
&\geq \frac{\nu_m}{\nu _{m+n}R^{m+n}}\int\limits_{B^n(b,R)}(R^2-\mid y-b\mid ^2)^{\frac{m}{2}}u_M(a,y)dm_n(y)\geq u_M(a,b).\end{align*}
Above we have used, in addition to the fact that, for every $y\in {\mathbb{R}}^m$,  the functions $u(\cdot ,y)$  are $K$-quasi-nearly
subharmonic, also the following lemma. (Observe that the proof of the Lemma, see [He71, proof of Theorem~2~a), p.~15], works also in our slightly more general situation: Recall that in the definition of nearly subharmonic functions, we use instead of the  standard condition $v\in {\mathcal{L}}^1_{\textrm{loc}}(U)$,  the slightly weaker condition $u^+\in {\mathcal{L}}^1_{\textrm{loc}}(U)$, see [Ri07$_4$, p.~51].)

\vspace{0.3ex}

\noindent{\textbf{Lemma.}} ([He71, Theorem~2~a), p.~15])  \emph{Let $v$ be  nearly subharmonic (in the generalized sense, defined above) in
a domain $U$ of ${\mathbb{R}}^{N}$, $N\geq 2$, $\psi \in {\mathcal{L}}^{\infty }({\mathbb{R}}^N)$, $\psi \geq 0$, $\psi (x)=0$ when $\mid x\mid \geq \alpha $ and $\psi (x)$ depends only on $\mid x\mid $.
 Then $\psi \star v\geq v$ and $\psi \star v$ is subharmonic in $U_{\alpha }$, provided $\int \psi (x)dm_N(x)=1$, where
$U_{\alpha }=\{x\in U:\, \overline{B^N(x,\alpha )}\subset U\}$.}
\hfill \qed
\section{An improvement to the result of Ko\l odziej and Thorbi\"ornson} 
\subsection{ } In our  generalization to the cited result of Ko\l odziej and Thorbi$\ddot {\textrm{o}}$rnson, Theorem~A above, we will use the generalized
Laplacian, defined with the aid of the Blaschke-Privalov operators, see e.g.  [Sa41, p.~451], [Ru50, pp.~278--279], [Sh56, p.~91],  [Br69, p.~20], [Sh71, p.~374],  [Sh78, p.~29] and [RiTa93, p.~1130]. Let  $D$ be a domain in
${\mathbb{R}}^{N}$,
 \mbox{$N\geq 2$,} and  $f:\, D \rightarrow {\mathbb{R}}$, $f\in {\mathcal{L}}^1_{\textrm{loc}}(D)$. We write
\begin{displaymath}\begin{split}
\Delta_*f(x):&=\liminf_{r\rightarrow 0}\frac{2(N+2)}{r^2}\cdot\big[\frac{1}{\nu _Nr^N}\int\limits_{B^N(x,r)}f(x')dm_N(x')-f(x)\big],\\
\Delta^* f(x):&=\limsup_{r\rightarrow 0}\frac{2(N+2)}{r^2}\cdot\big[\frac{1}{\nu _Nr^N}\int\limits_{B^N(x,r)}f(x')
dm_N(x')-f(x)\big].
\end{split}\end{displaymath}
If $\Delta _*f(x)= \Delta ^*f(x),$ then  write $\Delta f(x):= \Delta _*f(x)=\Delta ^*f(x)$.

If $f\in {\mathcal{C}}^2(D)$,
then
\begin{displaymath}\Delta f(x)=(\sum\limits_{j=1}^{N}\frac{\partial^2f}{\partial x_j^2})(x),
\end{displaymath}
the standard Laplacian with respect to the variable $x=(x_1,x_2,\dots ,x_N)$.
More generally, if $x\in D$ and
$f\in t^1_2(x)$, i.e. $f$ has an ${\mathcal{L}}^1$ total differential  of order $2$ at $x$,
then $\Delta f(x)$ equals with the pointwise Laplacian of $f$ at $x$, i.e.
\begin{displaymath}\Delta f(x)=\sum\limits_{j=1}^{N}D_{jj}f(x).\end{displaymath}
Here $D_{jj}f$ represents a generalization to the usual $\frac{\partial^2f}{\partial x_j^2}$,  $j=1,2, \dots ,N$.  See e.g. [CaZy61, p.~172], [Sh56, p.~498],
[Sh71, p.~369] and [Sh78, p.~29].

Recall that there exist functions which are not ${\mathcal{C}}^2$ but
for which the generalized Laplacian is nevertheless continuous, perhaps in the extended sense (in $([0,+\infty ],q)$, where $q$ is the spherical metric), see e.g. [Sh78, p.~31] and [Ri07$_4$, Example~5 and Example~6, pp.~66--67] and [Ri09, Example~1 and Example~2, pp.~e2624--e2625].

\vspace{0.3ex}

If $f$ is subharmonic on $D$, it follows from [Sa41, p.~451] (see also [Ru50, Lemma~2.2, p.~280]) that
$\Delta f(x)=\Delta_{*}f(x)=\Delta^*f(x)\in {\mathbb{R}}$
for almost every $x\in D$.

Below we use the following notation. Let  $\Omega $ is a domain in ${\mathbb{R}}^{m+n}$, \mbox{$m,n\geq 2$}, and
$u:\, \Omega \rightarrow {\mathbb{R}}$. If  $y\in {\mathbb{R}}^n$ is such that the function
\begin{displaymath} \Omega (y)\ni x\mapsto f(x):=u(x,y)\in {\mathbb{R}}\end{displaymath}
is in ${\mathcal{L}}^1_{\textrm{loc}}(\Omega (y))$, then we write $\Delta _{1*}u(x,y):=\Delta_*f(x)$,
$\Delta^* _1u(x,y):=\Delta^*f(x)$, and $\Delta _1u(x,y):=\Delta f(x)$.

\vspace{0.3ex}

\subsection{ } Then to our generalization to our previous result [Ri07$_4$, Theorem~5.1, pp.~67--68] (or  [Ri09, Theorem~4.3.1, p.~e2625] (where no proofs are given!)) and thus also to the result of  Ko\l odziej and Thorbi$\ddot {\textrm{o}}$rnson [KoTh96, Theorem~1, p.~463], Theorem~A above.   Though our proof will follow the main lines of  [Ri07$_4$, proof of Theorem~5.1, pp.~67--72], it is different enough, nevertheless, to warrant that it be given in complete detail here: Now  our  assumption (d) is  essentially milder than our  previous assumptions (d) and (e). 

\vspace{0.3ex}

\noindent{\textbf{Theorem~2}}  \emph{ Let $\Omega $ be a domain in ${\mathbb{R}}^{m+n}$, \mbox{$m,n\geq 2$.}
Let $u:\, \Omega \rightarrow
{\mathbb{R}}$ be such that for each $(x',y')\in \Omega $ there is  $(x_0,y_0)\in \Omega $ and $r_1>0$, $r_2>0$  such that
$(x',y')\in B^m(x_0,r_1)\times B^n(y_0,r_2)\subset \overline{B^m(x_0,r_1)\times B^n(y_0,r_2)}\subset \Omega $ and  such that the following
conditions are satisfied:}
 \begin{itemize}
\item[(a)] \emph{For each $y\in \overline{B^n(y_0,r_2)}$ the function}
\[\overline{B^m(x_0,r_1)}\ni x\mapsto u(x,y)\in {\mathbb{R}}\]
\emph{is continuous, and subharmonic in $B^m(x_0,r_1)$.}
 \item[(b)] \emph{For each $x\in \overline{B^m(x_0,r_1)}$ the function}
\[\overline{B^n(y_0,r_2)}\ni y\mapsto u(x,y)\in {\mathbb{R}}\]
\emph{is continuous, and harmonic in $B^n(y_0,r_2)$.}
\item[(c)] \emph{For each $y\in B^n(y_0,r_2)$ one has $\Delta _{1*}u(x,y)<+\infty $ for each $x\in B^m(x_0,r_1)$, possibly with the exception of
a polar set in $B^m(x_0,r_1)$.}
\item[(d)] \emph{There are a set $H\subset B^n(y_0,r_2)$, dense in $B^n(y_0,r_2)$, and a set $K\subset B^m(x_0,r_1)$, dense in $B^m(x_0,r_1)$, such that}
\begin{itemize}
\item[{(d1)}] \emph{for each $y\in H$, for almost every $x\in B^m(x_0,r_1)$ and for each $x\in K$,}
\[\Delta_{1}u(x',y)\rightarrow \Delta_1u(x,y)\in {\mathbb{R}}\]
\emph{as $x'\rightarrow x$, $x'\in K$, and}
\item[{(d2)}] \emph{for each $y\in B^n(y_0,r_2)\setminus H$ and for almost every $x\in B^m(x_0,r_1)$,} 
\[\Delta_{1}u(x,y')\rightarrow \Delta_1u(x,y)\in {\mathbb{R}}\]
\emph{as $y'\rightarrow y$, $y'\in H$.}
\end{itemize}
\end{itemize}
\emph{Then} $u$
\emph{is subharmonic in $\Omega $}.

\vspace{0.3ex}

\noindent\emph{Proof.} Choose $r'_1$,  $r'_2$ such that $0<r'_1<r_1$, $0<r'_2<r_2$, and  such that $(x',y')\in B^m(x_0,r'_1)\times B^n(y_0,r'_2)$.
It is sufficient to show that
$u\mid B^m(x_0,r'_1)\times B^n(y_0,r'_2)$ is subharmonic. For the sake of convenience of notation, we change the roles of $r_j$ and $r'_j$, $j=1,2$.
We divide the proof into several steps.

\vspace{0.3ex}

\noindent{\textbf{Step~1}} \emph{Construction of an auxiliar dense  set $G$.}

\vspace{0.3ex}

For each $k\in {\mathbb{N}}$ write
\[A_k:=\{\, x\in \overline{B^m(x_0,r_1)}\, :\, -k\leq u(x,y)\leq k \quad {\textrm{for each}}\quad y\in  \overline{B^n(y_0,r_2)}\,\}.\]
Clearly $A_k$ is closed, and
\[\overline{B^m(x_0,r_1)}=\bigcup _{k=1}^{+\infty }A_k.\]
Write
\begin{displaymath}G:=\bigcup _{k=1}^{+\infty }{\textrm{int}}\,\, A_k.\end{displaymath}
It follows from Baire's theorem  that $G$ is dense in $B^m(x_0,r_1)$.

\vspace{0.3ex}

\noindent{\textbf{Step~2}} \emph{The functions $\Delta_{1r}u(x,\cdot ): B^n(y_0,r_2)\rightarrow {\mathbb{R}}$ (see the definition below)}, $x\in G$,
$0<r<r_x:={\textrm{dist}}(x,\overline{B^m(x_0,r_1)}\setminus G)$, \emph{are
nonnegative and harmonic.}

\vspace{0.3ex}

For each $(x,y)\in B^m(x_0,r_1)\times B^n(y_0,r_2)$ and  each $r$, $0<r<{\textrm{dist}}(x,\partial B^m(x_0,r'_1))$ (observe that
${\textrm{dist}}(x,\partial B^m(x_0,r'_1))>r'_1-r_1>0$), write
\begin{displaymath}\begin{split}\Delta _{1r}u(x,y)&:=\frac{2(m+2)}{r^2}\cdot \big[\frac{1}{\nu _m\, r^m}\int\limits_{B^m(x,r)}u(x',y)\, dm_m(x')-u(x,y)\big]\\
&=\frac{2(m+2)}{r^2}\cdot \frac{1}{\nu _m\, r^m}\int\limits_{B^m(0,r)}\big[u(x+x',y)-u(x,y)\big]\, dm_m(x').\end{split}\end{displaymath}
Since $u(\cdot ,y)$ is subharmonic, $\Delta _{1r}u(x,y)$ is defined and nonnegative. Suppose then that $x\in G$ and $0<r<r_x$.
Since $\overline{B^m(x,r)}\subset G$ and $A_k\subset A_{k+1}$ for all $k=1,2, \dots$, $\overline{B^m(x,r)}\subset {\textrm{int}}\,\,A_N$ for some
$N\in {\mathbb{N}}$. Therefore
\begin{displaymath}-N\leq u(x',y)\leq N\quad {\textrm{for all}} \quad x'\in B^m(x,r) \quad {\textrm{and}} \quad y\in B^n(y_0,r_2),\end{displaymath}
and hence
\begin{equation}-2N\leq u(x+x',y)-u(x,y)\leq 2N\quad {\textrm{for all}} \quad x'\in B^m(0,r)\quad  {\textrm{and}} \quad y\in B^n(y_0,r_2).\end{equation}
To show that $\Delta_{1r}u(x,\cdot )$ is continuous, pick an arbitrary sequence $y_j\rightarrow \tilde{y}_0$, $y_j, \tilde{y}_0\in B^n(y_0,r_2)$,
$j=1,2, \dots$. Using then (1), Lebesgue Dominated Convergence Theorem and the continuity of $u(x, \cdot )$, one sees easily that
$\Delta _{1r}u(x,\cdot )$ is continuous.

To show that $\Delta _{1r}u(x,\cdot )$ satisfies the mean value
equality, take $\tilde{y}_0\in B^n(y_0,r_2)$ and $\rho >0$ arbitrarily such
that $\overline{B^n(\tilde{y}_0,\rho )}\subset B^n(y_0,r_2)$.
Because of (1) we can use Fubini's Theorem. Thus
\begin{displaymath}\begin{split}\frac{1}{\nu _n\rho ^n}\int\limits_{B^n(\tilde{y}_0,\rho )}\Delta _{1r}u(x,y) dm_n(y)
&=\frac{1}{\nu _n\rho ^n}\int\limits_{B^n(\tilde{y}_0,\rho )}\{\frac{2(m+2)}{r^2}\cdot \frac{1}{\nu _m
 r^m}\int\limits_{B^m(0,r)}\big[u(x+x',y)-u(x,y)\big] dm_m(x')\}dm_n(y)\\
&=\frac{2(m+2)}{r^2}\cdot \frac{1}{\nu _m
 r^m}\int\limits_{B^m(0,r)}\{\frac{1}{\nu _n\rho ^n}\int\limits_{B^n(\tilde{y}_0,\rho )}\big[u(x+x',y)-u(x,y)\big]dm_n(y)\} dm_m(x')\\
&=\frac{2(m+2)}{r^2}\cdot \frac{1}{\nu _m\, r^m}\int\limits_{B^m(0,r)}\big[u(x+x',\tilde{y}_0)-u(x,\tilde{y}_0)\big]\, dm_m(x')\\
&=\Delta _{1r}u(x,\tilde{y}_0).\end{split}\end{displaymath}

\vspace{0.3ex}

\noindent{\textbf{Step~3}} \emph{The functions  $\Delta_{1}u(x,\cdot ): B^n(y_0,r_2)\rightarrow {\mathbb{R}}$, $x\in G\cap A$,  are defined, nonnegative and harmonic. Here}
\begin{equation*}A:=\bigcap_{k=1}^{+\infty}A(y_k),\end{equation*}
\emph{where $H=\{\, y_k, \, k=1,2, \dots \, \}$ (we may clearly suppose that $H$ is countable), and, for arbitrary $y\in B^n(y_0,r_2)$,}
\[A(y):=\{\,x\in B^m(x_0,r_1):\,\Delta_{1*}u(x,y)= \Delta_1^*u(x,y) =\Delta_1u(x,y)\in  {\mathbb{R}}\,\}.\]

\vspace{0.3ex}

By [Ru50, Lemma~2.2, p.~280] (see also [Sa41, p.~451] and [Sh71, p.~376]), $m_m(B^m(x_0,r_1)\setminus A(y))=0$ for each $y\in B^n(y_0,r_2)$.

Take $x\in G\cap A$ and a sequence  $r_j\rightarrow 0$, $0<r_j<r_x$, $j=1,2,\dots$, arbitrarily. By
[He71, Corollary~3 a), p.~6] (or [AG01, Lemma~1.5.6, p.~16]) we see that the family
\begin{displaymath} \Delta_ {1r_j}u(x,\cdot ):\, B^n(y_0,r_2) \rightarrow {\mathbb{R}}, \, j=1,2,\dots ,\end{displaymath}
of nonnegative and harmonic functions is either uniformly equicontinuous and locally  uniformly bounded, or else
\[\sup_{j=1,2,\dots}\Delta_{1r_j}u(x,\cdot )\equiv +\infty .\]
On the other hand, since $x\in G\cap A$, we know that for each $y_k\in H, \,k=1,2,\dots$, 
\begin{displaymath}\Delta_{1r_j}u(x,y_k)\rightarrow \Delta_1u(x,y_k)\in {\mathbb{R}}\, 
\end{displaymath}
as $j\rightarrow +\infty$.
Therefore,  by [V\"a71, Theorem~20.3, p.~68] and by [He71, c), p.~2] (or [ArGa01,  Theorem~1.5.8, p.~17]), the limit 
\begin{displaymath}\Delta_1u(x,\cdot )=\lim_{j\rightarrow +\infty }\Delta_{1r_j}u(x,\cdot )\end{displaymath}
exists and defines a harmonic function in $B^n(y_0,r_2)$. Since the limit is clearly independent of the considered sequence $r_j$, the claim follows.

\vspace{0.3ex}

\noindent{\textbf{Step~4}} \emph{The function  $\Delta_{1}u(\cdot
,\cdot )\mid (G\cap K \cap A\cap B)\times B^n(y_0,r_2)$ has a continuous  extension
$\tilde{\Delta}_{1}u(\cdot ,\cdot ): (A\cap B)\times
B^n(y_0,r_2)\rightarrow {\mathbb{R}}$.   Moreover, the functions
$\tilde{\Delta}_{1}u(x,\cdot ): 
B^n(y_0,r_2)\rightarrow {\mathbb{R}}$, $x\in A\cap  B$, are nonnegative and
harmonic. Here}
\begin{equation*}B:=\bigcap_{k=1}^{+\infty}B(y_k),\end{equation*}
\emph{where, for arbitrary $y\in B^n(y_0,r_2)$,  we use the notation}
\[B(y):=\{\,x\in B^m(x_0,r_1):\,\Delta_{1}u(x',y) \rightarrow \Delta_1u(x,y) {\textrm{ as }} x'\rightarrow x, \,\, x'\in K\,\}.\]

\vspace{0.3ex}

Using the assumption (d1), one sees easily that $G\cap K\cap A\cap B$ is dense in $A\cap B$.

To show the existence of the desired continuous extension, it is clearly sufficient to show that for each $(\tilde{x}_0,\tilde{y}_0)\in
(A\cap B)\times B^n(y_0,r_2)$, the limit
\begin{displaymath}\lim_{(x,y)\rightarrow (\tilde{x}_0,\tilde{y}_0),\, (x,y)\in (G\cap K\cap A\cap B)\times B^n(y_0,r_2)} \Delta _1u(x,y)\end{displaymath}
exists.  (This is of course standard, see e.g. [Di60, (3.15.5), p.~54].) To see this, it is  sufficient to show that,
for an arbitrary sequence $(x_j,y_j)\rightarrow (\tilde{x}_0,\tilde{y}_0)$, $(x_j,y_j)\in (G\cap K\cap A\cap B)\times B^n(y_0,r_2)$, $j=1,2, \dots$,  the limit
\begin{displaymath}\lim_{j\rightarrow +\infty }\Delta _1u(x_j,y_j)\end{displaymath}
exists.

That this limit indeed exists, is seen as above, just using the facts:
\begin{itemize}
\item[${-}$] the functions $\Delta _1u(x_j, \cdot )$, $j=1,2,\dots$, are nonnegative and harmonic in $B^n(y_0,r_2)$, by Step~3;
\item[${-}$] for each $y_k\in H$, $k=1,2,\dots$,   $\Delta _1u(x_j,y_k)\rightarrow \Delta_{1}u(\tilde{x}_0,y_k)\in {\mathbb{R}}$  as $j\rightarrow +\infty $.
\end{itemize}
(See again [He71, Corollary~3 a), p.~6] (or [AG01, Lemma~1.5.6, p.~16]) and [V\" a71, Theorem~20.3,  p.~68]).
That the functions $\tilde{\Delta}_1u(x,\cdot ):\, B^n(y_0,r_2)\rightarrow {\mathbb{R}}$, $x\in A\cap B$, are harmonic, see [He71, c), p.~2] (or [ArGa01,  Theorem~1.5.8, p.~17])).
 
\vspace{0.3ex}

\noindent{\textbf{Step~5}} \emph{For each $x\in B^m(x_0,r_1)$ the  functions}
\begin{displaymath}B^n(y_0,r_2)\ni y\mapsto  \tilde{v}(x,y):=\int G_{B^m(x_0,r_1)}(x,z)\tilde{\Delta }_1u(z,y)dm_m(z)\in {\mathbb{R}}\end{displaymath}
\emph{and}
\begin{displaymath}B^n(y_0,r_2)\ni y\mapsto  \tilde{h}(x,y):=u(x,y)+\tilde{v}(x,y)\in {\mathbb{R}}\end{displaymath}
\emph{are harmonic. Above and below $G_{B^m(x_0,r_1)}(x,z)$ is the Green function of the ball $B^m(x_0,r_1)$, with $x$ as a pole.}

\vspace{0.3ex}

Using  Fubini's Theorem one sees easily that for each $x\in B^m(x_0,r_1)$ the function
 $\tilde{v}(x, \cdot )$ satisfies the mean value equality. To see that $\tilde{v}(x, \cdot )$ is harmonic, it is sufficient to show that
$\tilde{v}(x, \cdot )\in {\mathcal{L}}^1_{\textrm{loc}}(B^n(y_0,r_2))$.   Using just   Fatou's Lemma, one sees that $\tilde{v}(x, \cdot )$ is
lower semicontinuous, hence superharmonic.
Therefore either $\tilde{v}(x,\cdot )\equiv +\infty $ or else $\tilde{v}(x, \cdot )\in {\mathcal{L}}^1_{\textrm{loc}}(B^n(y_0,r_2))$. The following argument shows
that the former alternative cannot
occur. Indeed,   for each $x\in A\cap B$ and
for each $y_k\in H$, $k=1,2,\dots$, we see, using the definition of the (continuous) function $\tilde{\Delta}_1u(\cdot ,\cdot )$ and (d1), that 
\begin{equation}\tilde{\Delta}_{1}u(x,y_k)=\lim_{x'\rightarrow x,\,x'\in G\cap K\cap A\cap B}\tilde{\Delta}_{1}u(x',y_k)=\lim_{x'\rightarrow x,\,x'\in G\cap K\cap A\cap B}
\Delta_{1}u(x',y_k)=\Delta_{1}u(x,y_k)\in {\mathbb{R}}.\end{equation}
 Hence $\tilde{v}(x,y_k)=v(x,y_k)\in {\mathbb{R}}$  for each $x\in B^m(x_0,r_1)$ and $y_k\in H$, $k=1,2,\dots$. (See (3) in Step~6 below for the definition of $v(\cdot ,\cdot ): B^m(x_0,r_1)\times B^n(y_0,r_2) \rightarrow {\mathbb{R}}$.)  Therefore, for each $x\in B^m(x_0,r_1)$, the function $\tilde{v}(x,\cdot )$ and thus also  the function $\tilde{h}(x,\cdot )=u(x,\cdot )+\tilde v(x,\cdot )$ are harmonic.

\vspace{0.3ex}

\noindent{\textbf{Step~6}} \emph{For each $y\in B^n(y_0,r_2)$ the  function}
\begin{displaymath}B^m(x_0,r_1)\ni x\mapsto  \tilde{h}(x,y)\in {\mathbb{R}}\end{displaymath}
\emph{is harmonic.}

\vspace{0.4ex}

With the aid of the version of Riesz's Decomposition Theorem, given in  [Ru50, 1.3. Theorem~II, p.~279, and p.~278, too]
(see also [Sh56, Theorem~1, p.~499]), for  each $y\in B^n(y_0,r_2)$ one can write
\begin{displaymath}u(x,y)=h(x,y)-v(x,y),\end{displaymath}
where
\begin{equation}v(x,y):=\int G_{B^m(x_0,r_1)}(x,z)\Delta _1u(z,y)dm_m(z)\end{equation}
and $h(\cdot ,y)$ is the least harmonic majorant of $u(\cdot ,y)\mid B^m(x_0,r_1)$. Here  $v(\cdot ,y)$ is continuous and superharmonic in $B^m(x_0,r_1)$.

As shown above in  (2),  $v(\cdot ,y_k)=\tilde{v}(\cdot ,y_k)$ for each $y_k\in H$, $k=1,2,\dots $. Therefore $\tilde{h}(\cdot, y_k)=h(\cdot ,y_k)$, and thus $\tilde{h}(\cdot, y_k$) is harmonic for each $y_k\in H$, $k=1,2,\dots$.  

To see that $\tilde{h}(\cdot ,y)$ is harmonic also for  $y\in B^n(y_0,r_2)\setminus H$, take $\tilde{y}_0\in B^n(y_0,r_2)\setminus H$ arbitrarily, and 
proceed in the following way. Take  $z\in A\cap B\cap A(\tilde{y}_0)\cap C(\tilde{y}_0)$ arbitrarily, where, for arbitrary $y\in B^n(y_0,r_2)\setminus H$, 
\[C(y):=\{\,z\in B^m(x_0,r_1):\,\Delta_{1*}u(z,y') \rightarrow \Delta_{1*}u(z,y) {\textrm{ as }} y'\rightarrow y, \,\ y'\in H\, \}.\]
Since $z\in A(\tilde{y}_0)$, we have $\Delta_{1*}u(z,\tilde{y}_0)=\Delta_{1}u(z,\tilde{y}_0)\in {\mathbb{R}}$. Thus we may also suppose that $\Delta_{1*}u(z,y')=\Delta_{1}u(z,y')\in {\mathbb{R}}$. Using then our assumption (d2) and the continuity of  $\tilde{\Delta}_1u(\cdot ,\cdot )$, we see that  
\begin{displaymath}\Delta_1u(z,\tilde{y}_0)=\tilde{\Delta}_1u(z,\tilde{y}_0) \end{displaymath}
for every $z\in A\cap B\cap A(\tilde{y}_0)\cap C(\tilde{y}_0)$. Therefore,  $\tilde{v}(x,\tilde{y}_0)=v(x,\tilde{y}_0)$ and thus $\tilde{h}(x,\tilde{y}_0)=h(x,\tilde{y}_0)$ for each $x\in B^m(x_0,r_1)$. 

\vspace{0.3ex}

\noindent{\textbf{Step~7}} \emph{The use of the results of  Lelong and of Avanissian.}

\vspace{0.3ex}

By Steps~5 and 6 we know that $\tilde{h}(\cdot ,\cdot )=h(\cdot ,\cdot )$ is separately harmonic in $B^m(x_0,r_1)\times B^n(y_0,r_2)$.
By Lelong's result [Le61, p.~561] (see also [Av67, Théorème~1, pp.~4--5]) $\tilde{h}(\cdot ,\cdot )$ is harmonic and thus
locally bounded above in $B^m(x_0,r_1)\times B^n(y_0,r_2)$. Therefore also $u(\cdot ,\cdot )$ is locally bounded above in $B^m(x_0,r_1)\times B^n(y_0,r_2)$. 
But then it follows from  [Av61, Théorème~9, p. 140] (see also [Le45, Théorème~1 bis, p.~315],  [Ar66, Theorem~1, p.~622], [Le69, Proposition~3, p.~24],  [Ri89,Theorem~1, p. 69], [ArGa93, Theorem~1, p.~256] and [Ri08, Corollary~4.5, p.~13]) that $u(\cdot ,\cdot )$ is subharmonic on $B^m(x_0,r_1)\times B^n(y_0,r_2)$.
\hfill \qed
\vspace{0.3ex}

\noindent{\textbf{Remark~1}} Observe that the assumption (d2) was needed \emph{only} to see that 
\[\Delta_1 u(x,y) = \tilde{\Delta}_1 u(x,y) {\textrm{ for almost every }} x\in B^m(x_0,r_1){\textrm{ and for each }} y\in B^n(y_0,r_2)\setminus H .\] 
(At this point one might recall that the functions $\tilde{\Delta}_1u(x,\cdot ):\, B^n(y_0,r_2) \rightarrow {\mathbb{R}}$, $x\in B^m(x_0,r_1)$, are harmonic.)

From the above proof one sees easily that the  assumption (d), that is (d1) and (d2), can be replaced by:
\begin{itemize}
\item[{(d$^*$)}]  \emph{For every $y\in B^n(y_0,r_2)$, for almost every $x\in B^m(x_0,r_1)$ and for each $x\in K$,} 
\[\Delta_{1}u(x',y)\rightarrow \Delta_1u(x,y)\in {\mathbb{R}}\]
\emph{as $x'\rightarrow x$, $x'\in K$.}
\end{itemize}

\vspace{0.3ex}

Though our Theorem~2 might still be considered somewhat technical, it has, nevertheless, the following concise  corollaries, both of which already  improve the result of Ko\l odziej and J.~Thorbi$\ddot {\textrm{o}}$rnson. 

\vspace{0.3ex}

\noindent{\textbf{Corollary~1}} ([Ri07$_4$, Corollary~5.1, p.~74] and  [Ri09, Corollary~4.3.3,  p.~e2626]) \emph{ Let $\Omega $ be a domain in ${\mathbb{R}}^{m+n}$, \mbox{$m,n\geq 2$.}
Let $u:\, \Omega \rightarrow
{\mathbb{R}}$ be such that}
 \begin{itemize}
\item[(a)] \emph{for each $y\in {\mathbb{R}}^n$ the function}
\[\Omega (y)\ni x\mapsto u(x,y)\in {\mathbb{R}}\]
\emph{is continuous and subharmonic,}
 \item[(b)] \emph{for each $x\in {\mathbb{R}}^m$ the function}
\[\Omega (x)\ni y\mapsto u(x,y)\in {\mathbb{R}}\]
\emph{is harmonic,}
\item[(c)] \emph{for each $y\in {\mathbb{R}}^n$ the function}
\[\Omega (y)\ni x\mapsto \Delta _1u(x,y)\in [0,+\infty ]\]
\emph{is defined,  continuous (with respect to the spherical metric), and finite for all $x$, except at most of a polar set $E(y)$ in $\Omega (y)$.}
\end{itemize}
\emph{Then} $u$
\emph{is subharmonic in $\Omega $}.

\vspace{0.3ex}

\noindent{\textbf{Corollary~2}} ([Ri07$_3$, Corollary, p.~444])  \emph{ Let $\Omega $ be a domain in ${\mathbb{R}}^{m+n}$, \mbox{$m,n\geq 2$.}
Let $u:\, \Omega \rightarrow
{\mathbb{R}}$ be such that}
 \begin{itemize}
\item[(a)] \emph{for each $y\in {\mathbb{R}}^n$ the function}
\[\Omega (y)\ni x\mapsto u(x,y)\in {\mathbb{R}}\]
\emph{is continuous and subharmonic,}
 \item[(b)] \emph{for each $x\in {\mathbb{R}}^m$ the function}
\[\Omega (x)\ni y\mapsto u(x,y)\in {\mathbb{R}}\]
\emph{is harmonic,}
\item[(c)] \emph{for each  $y\in {\mathbb{R}}^n$ the function}
\[\Omega (y)\ni x\mapsto \Delta _1u(x,y)\in {\mathbb{R}}\]
\emph{is defined  and continuous.}
\end{itemize}
\emph{Then} $u$
\emph{is subharmonic in $\Omega $}.

\vspace{0.3ex}

\flushleft{\noindent\textbf{References}}

\vspace{0.2ex}

\begin{flushleft}
\begin{enumerate}
\item[{[ArGa93]}] D.H.~Armitage and S.J.~Gardiner,
 \emph{Conditions for separately subharmonic functions to be subharmonic},
Potential Anal.,
{\textbf{2}} (1993), 255--261.
\item[{[ArGa01]}] D.H.~Armitage and S.J.~Gardiner,
 \emph{Classical Potential Theory},
Springer-Verlag, London, 2001.
\item[{[Ar66]}] M.G.~Arsove,  \emph{On subharmonicity of doubly subharmonic functions},
Proc. Amer. Math. Soc.,
{\textbf{17}} (1966), 622--626.
\item[{[Av61]}] V.~Avanissian,
 \emph{Fonctions plurisousharmoniques et fonctions doublement sousharmoniques},
Ann. Sci. École  Norm. Sup., {\textbf{78}} (1961), 101--161.
\item[{[Av67]}] V.~Avanissian, 
 \emph{Sur l'harmonicité des fonctions séparément harmoniques}, in: 
Séminaire de Probabilités (Univ. Strasbourg, Février 1967), {\textbf{1}} (1966/1967), pp.~101--161, Springer, Berlin, 1967.
\item[{[Br69]}] M.~Brelot,
 \emph{Éléments de la Théorie Classique du Potentiel},
Centre de Documentation Universitaire, Paris, 1969 (Third Edition).
\item[{[CaZy61]}] A.P.~Calderon and A.~Zygmund,  \emph{Local properties of solutions of elliptic partial differential equations},
Studia Math., {\textbf{20}} (1961), 171--225.
\item[{[CeSa93]}] U.~Cegrell and A.~Sadullaev, \emph{Separately subharmonic functions}, Uzbek. Math. J., {\textbf{1}} (1993), 78--83.
\item[{[Di60]}] J.~Dieudonné,
 \emph{Foundations of Modern Analysis},
Academic Press, New York, 1960.
\item[{[He71]}] M.~Hervé,  \emph{Analytic and Plurisubharmonic Functions in Finite and Infinite Dimensional Spaces},
Lecture Notes in Mathematics, Vol. {\textbf{198}}, Springer, Berlin $\cdot$ Heidelberg $\cdot$ New York, 1971.
\item[{[Im90]}] S.A.~Imomkulov, \emph{Separately subharmonic functions} (in Russian),
Dokl. USSR, {\textbf{2}} (1990), 8--10.
\item[{[KoTh96]}] S.~Ko\l odziej and J.~Thorbi$\ddot {\textrm{o}}$rnson,  \emph{Separately harmonic and subharmonic functions},
Potential Anal., {\textbf{5}} (1996), 463--466.
\item[{[Le45]}] P.~Lelong, \emph{Les fonctions plurisousharmoniques},
Ann. Sci. École Norm. Sup.,
{\textbf{62}} (1945), 301--338.
\item[{[Le61]}] P.~Lelong, \emph{Fonctions plurisousharmoniques et fonctions analytiques de variables réelles},
Ann. Inst. Fourier, Grenoble,
{\textbf{11}} (1961), 515--562.
\item[{[Le69]}] P.~Lelong, \emph{Plurisubharmonic Functions and Positive Differential Forms},
Gordon and Breach, London, 1969.
\item[{[PaRi08]}] M.~Pavlovi\'c  and J.~Riihentaus,  \emph{Classes of quasi-nearly subharmonic functions}, Potential Anal. {\textbf{29}} (2008), 89--104.
\item[{[PaRi09]}] M.~Pavlovi\'c and J.~Riihentaus, \emph{Quasi-nearly subharmonic functions in locally uniformly homogeneous spaces},  Positivity, \textbf{15}, no.~{\textbf{1}} (2009), 1-10. 
\item[{[Ri89]}] J.~Riihentaus, \emph{On a theorem of Avanissian--Arsove},
Expo.  Math., {\textbf{7}} (1989),  \mbox{69--72.}
\item[{[Ri07$_1$]}] J.~Riihentaus, \emph{Separately quasi-nearly subharmonic functions}, in: Complex Analysis and Potential Theory, Proceedings of the
Conference Satellite to ICM~2006, Tahir Aliyev Azero$\breve{\textrm{g}}$lu, Promarz M. Tamrazov (eds.), Gebze Institute of Technology, Gebze, Turkey,
September 8-14,  2006, World Scientific, Singapore, 2007, pp.~156--165.
\item[{[Ri07$_2$]}] J.~Riihentaus, \emph{On the subharmonicity of separately  subharmonic functions}, in: Proceedings of the 11th WSEAS International
Conference on Applied Mathematics (MATH'07),  Dallas, Texas, USA, March 22-24, 2007, Kleanthis Psarris, Andrew D.~Jones (eds.), WSEAS, 2007,
pp.~230-236. 
\item[{[Ri07$_3$]}] J.~Riihentaus, \emph{On separately harmonic and subharmonic functions}, Int. J. Pure Appl. Math., {\textbf{35}}, no. {\textbf{4}} (2007),
\mbox{435-446}.
\item[{[Ri07$_4$]}] J.~Riihentaus, \emph{Subharmonic functions, generalizations and separately subharmonic functions},
The XIV-th Conference on Analytic Functions, July 22-28, 2007, Che\l m, Poland, in: Scientific Bulletin of Che\l m, Section of Mathematics
and Computer Science, {\textbf{2}} (2007), 49--76.  
 \item[{[Ri08]}] J.~Riihentaus,  \emph{Quasi-nearly subharmonicity and separately quasi-nearly subharmonic functions}, J. Inequal. Appl.,
{\textbf{2008}}, Article ID 149712, 15~pages, 2008. 
\item[{[Ri09]}] J.~Riihentaus, \emph{Subharmonic functions, generalizations and separately subharmonic functions: A survey},
5th World Congress of Nonlinear Analysts (WCNA '08), 
July 2 - 9, 2008, Orlando, Florida, USA, in: Nonlinear Analysis, {\textbf{71}} (2009), e2613--e2627. 
\item[{[Ri11]}] J.~Riihentaus, \emph{Domination conditions for families of quasinearly subharmonic functions},  International Journal of Mathematics and Mathematical Sciences/\textit{New Trends in Geometric Function Theory~2011}, {\textbf{2011}} (2011), Article ID~729849, 9~pages.
\item[{[RiTa93]}] J.~Riihentaus and P.M.~Tamrazov,  \emph{On subharmonic extension and the extension in the Hardy-Orlicz classes} (English and Ukrainian summaries), 
Ukrain. Mat. Zh., {\textbf{45}}, no. {\textbf{8}} (1993), 1260--1271. 
\item[{[Ru50]}] W.~Rudin, \emph{Integral representation of continuous functions}, Trans. Amer. Math. Soc., {\textbf{68}} (1950),
278--286.
\item[{[Sa41]}] S.~Saks, \emph{On the operators of Blaschke and Privaloff for subharmonic functions}, Rec. Math. (Mat. Sbornik), {\textbf{9 (51)}} (1941),
451--456.
\item[{[Sh56]}] V.L.~Shapiro, \emph{Generalized laplacians},  Amer. J. Math.,
{\textbf{78}} (1956), \mbox{497--508.}
\item[{[Sh71]}] V.L.~Shapiro, \emph{Removable sets for pointwise subharmonic functions}, Trans. Amer. Math. Soc.,
{\textbf{159}} (1971), \mbox{369--380.}
\item[{[Sh78]}] V.L.~Shapiro, \emph{Subharmonic functions and Hausdorff measure}, J. Diff. Eq., {\textbf{27}} (1978), \mbox{28--45.}
\item[{[Ta88]}] P.M.~Tamrazov,  \emph{Removal of singularities of subharmonic, plurisubharmonic functions and their generalizations} (English and Ukrainian summaries), 
Ukrain. Mat. Zh., {\textbf{40}}, no. {\textbf{6}} (1988), 683--694 (Russian); translation in Ukrainian Math. J., {\textbf{40}}, no. {\textbf{6}} (1988), 573--582. 
\item[{[V\" a71]}] J.~V\"ais\"al\"a, \emph{Lectures on n-Dimensional Quasiconformal Mappings},
Lecture Notes in Mathematics, Vol. {\textbf{229}}, Springer, Berlin $\cdot$ Heidelberg $\cdot$ New York, 1971.
\item[{[Wi88]}] J.~Wiegerinck, \emph{Separately  subharmonic functions need not be subharmonic},
Proc.  Amer. Math.  Soc., {\textbf{104}} (1988), 770--771.
\item[{[WiZe91]}] J.~Wiegerinck and R.~Zeinstra, \emph{Separately subharmonic functions: when are they subharmonic}, in: Proceedings of Symposia in Pure
Mathematics, vol. {\textbf{52}}, part {\textbf{1}},  Eric Bedford, John P. D'Angelo, Robert E.~Greene, Steven G.~Krantz (eds.),
 Amer. Math. Soc., Providence, Rhode Island, 1991,  \mbox{pp. 245--249}.
\end{enumerate}
\end{flushleft}
\end{document}